\documentclass[reprint, amsmath,amssymb, aps,
 prl]{revtex4-2}

\usepackage{graphicx}% Include figure files
\usepackage{dcolumn}% Align table columns on decimal point
\usepackage{bm}% bold math

\usepackage{xcolor}

\newcommand{\tp}[1]{#1^\mathsf{T}}

\begin{document}

\preprint{APS/123-QED}

\title{Phase and amplitude responses for delay equations using harmonic balance}

\author{R. Nicks}
\email{Rachel.Nicks@nottingham.ac.uk}
\author{R. Allen}
\email{pcyrga@exmail.nottingham.ac.uk}
\author{S. Coombes}
\email{Stephen.Coombes@nottingham.ac.uk}
%\homepage[]{Your web page}
%\thanks{}
%\altaffiliation{}
\affiliation{School of Mathematical Sciences, University of Nottingham, Nottingham, NG7 2RD, UK}

\date{\today}% It is always \today, today,
             %  but any date may be explicitly specified

\begin{abstract}
Robust delay induced oscillations, common in nature, are often modeled by delay-differential equations (DDEs).  Motivated by the success of phase-amplitude reductions for ordinary differential equations with limit cycle oscillations, there is now a growing interest in the development of analogous approaches for DDEs to understand their response to external forcing.  When combined with Floquet theory, the fundamental quantities for this reduction are phase and amplitude response functions.  Here, we develop a framework for their construction that utilises the method of harmonic balance.  
%593 char
\end{abstract}

%\keywords{Suggested keywords}%Use showkeys class option if keyword
                              %display desired
\maketitle

%\tableofcontents

Time delays can lead to oscillatory behaviour in many real world systems, as exemplified by laser networks \cite{Erneux2009}, machine dynamics \cite{Stepan1989}, and neural systems \cite{Campbell2007}.
These can often be described by delay-differential equations (DDEs) with a stable limit cycle. Such systems do not occur in isolation and are perturbed by external forcing, or interactions with other oscillators. It is therefore vital to be able to quantify the effect these perturbations have on the system behaviour.  However, the understanding of the response of DDEs to external forces is challenging due to their infinite dimensionality

For weak perturbations a popular method of oscillator reduction for ordinary differential equations (ODEs) is based on formulating phase dynamics along a cycle, and this has recently been extended to account for some notion of distance from cycle using isostable coordinates \cite{Guillamon2009,Mauroy2013,Wilson2016}.  For a recent overview see \cite{Ermentrout2019}.  
It has been shown that phase-amplitude reduction retaining the phase and slowest decaying amplitude (isostable) of oscillators can capture qualitative changes in stability of phase-locked states occurring under increasing interaction strength in networks of coupled ODE oscillators \cite{Nicks2024}. 
Extending these results to networks of delayed systems requires the computation of phase and amplitude response to perturbations of limit cycle solutions of DDEs, as well as associated Floquet exponents and eigenfunctions. 
The phase only reduction has previously been generalised to DDEs in \cite{Novicenko2012,Kotani2012}, and only recently has a phase-amplitude formulation been proposed using a functional analytic perspective by Kotani \textit{et al}. \cite{Kotani2020}.
In both these approaches the practical application of the theory is via the numerical solution of DDEs using time evolution methods.  Here, we propose an alternative approach, based upon the method of harmonic balance, that can side-step the numerical challenges associated with time evolution.

In this paper we describe how to arrive at the linear (adjoint) equations for phase and amplitude responses (with appropriate normalisations) using a generalisation of the approach of Novi\v{c}enko and Pyragas \cite{Novicenko2012}. This involves a discretisation of the (infinite dimensional) DDE system as a system of (high dimensional) ODEs and the use of the ODE theory for phase-amplitude reduction to obtain the equations for the DDE phase and amplitude responses in the continuum limit. We further introduce a practical methodology, based upon harmonic balance, to approximate the periodic solutions of the (linear) DDEs describing the response functions in addition to determining the delay induced orbit (from a nonlinear DDE) and its Floquet exponents and eigenfunctions.  To illustrate the utility of our approach we compare against two nonlinear DDE models for which analytical results are known.

For limit-cycle oscillators described by ODEs of the form
\begin{equation}\label{eq:ODE} \dot{y} = G(y)+ \epsilon P(t),\end{equation} where $y \in \mathbb{R}^m$ and $\epsilon P(t)$ is a small time-dependent perturbation, we assume that when $\epsilon=0$ the system has a stable limit cycle $y^\gamma(t)$ with period $T$. The phase $\theta = \theta(y)$ and amplitudes $\psi_i= \psi_i(y)$ can be defined near the limit cycle such that $\dot{\theta}=\omega = 2\pi/T$ and $\dot{\psi}_i = \mu_i \psi_i$ where $\mu_i$, $i=1, \cdots m-1$ are Floquet exponents for which $\mbox{Re}(\mu_i)<0$. The exponent with real part closest to zero is $\mu_1:=\mu$ which here we assume to be real and small so that perturbations in the direction of the corresponding eigenfunction $g_1(t):=g(t)$ decay slowly. We assume that all other exponents have $\mbox{Re}(\mu_i)$ large and negative so we may retain only $\psi_1 :=\psi$ and all other amplitudes may be neglected as they decay much faster.

In the presence of forcing, the standard first order phase-amplitude reduction is given by \cite{Wilson2018},
\begin{align}
\dot{\theta} &= \omega + \epsilon \tp{Z}(t) P(t) ,\\
\dot{\psi} &= \mu \psi + \epsilon \tp{I}(t) P(t) ,
\end{align} where $\tp{}$ denotes transpose, but higher order corrections may also be included \cite{Wilson2018, Park2021, Nicks2024}. Here $Z(t)$ and $I(t)$ are the $T$-periodic phase and amplitude response functions respectively, evaluated on the limit cycle. They quantify the linear response of the phase and amplitude of the oscillator to the perturbation and  can be computed as the $T$ periodic solutions of the adjoint equations
\begin{align}
\dot{Z}(t) &= -\tp{J} Z(t) ,\label{eq:adjointZ}\\
\dot{I}(t) &= -(\tp{J} - \mu I_m) I(t) ,\label{eq:adjointI}
\end{align} normalised according to 
\begin{equation}
\tp{Z}(0)\dot{y}^\gamma(0)=\omega , \qquad
\tp{I}(0)g(0) =1,
\end{equation} where $I_m$ is the $m \times m$ identity matrix and $J:={\rm D}G(y^\gamma(t))$ is the Jacobian of the vector field $G$ evaluated on the limit cycle. Here $g(t)$ is the $T$-periodic Floquet eigenfunction associated with the exponent $\mu$ which satisfies the linear equation \cite{Wilson2020b}
\begin{align}\label{eq:FloquetODE}
\dot{g}(t) = (J - \mu I_m) g(t) ,
\end{align} and we are free to specify the normalisation, usually making $|g(0)|=1$. 

In the current work we focus on the computation of phase and amplitude response functions for perturbed limit cycle oscillators described by DDEs of the form
\begin{equation}
\dot{x}(t) = F(x(t), x(t-\tau)) + \epsilon p(t) , \label{eq:DDE}
\end{equation} where $x\in\mathbb{R}^m$ and $\tau$ is a constant delay. We assume that when $\epsilon=0$, \eqref{eq:DDE} admits a limit cycle solution $x^{\gamma}(t)$ with period $T$. Following \cite{Novicenko2012}, equation \eqref{eq:DDE} is equivalent to 
\begin{align}
\dot{x}(t) &= F(x(t), \xi(\tau, t)) + \epsilon p(t)\\
\frac{\partial \xi(s,t)}{\partial t} &= -\frac{\partial \xi(s,t)}{\partial s}, \quad \xi(0,t) = x(t), 
\end{align} where $s \in [0, \tau]$ and therefore $\xi(\tau, t)= x(t-\tau)$.
Next, discretise this as a system of $m(N+1)$ ODEs by defining $x_0(t) = x(t)$ and $x_i(t) = \xi(i\tau/N,t)\approx x(t-i\tau/N)$ for $i=1,\ldots,N$. Then 
\begin{align}
\dot{y}(t) = G(y(t))+\epsilon P(t) , \label{eq:Discretised}
\end{align}
where $y(t) = \tp{(\tp{x_0}(t), \tp{x_1}(t), \ldots \tp{x_N}(t))}$,  $P(t) = \tp{(\tp{p}(t),0, \ldots, 0)}$,
\begin{align*}
G(y(t)) = 
\begin{pmatrix}
F(x_0(t),x_N(t))\\
N[x_0(t)-x_1(t)]/\tau \\
\vdots \\
N[x_{N-1}(t)-x_N(t)]/\tau
\end{pmatrix}.
\end{align*}
System \eqref{eq:Discretised} with $\epsilon=0$ has Jacobian
\begin{align} \label{eq:DiscretisedJ}
J(t) = \begin{pmatrix}
{\rm D}F_0(t) & 0 & 0 & \cdots & {\rm D}F_1(t) \\
\frac{N}{\tau} I_m & -\frac{N}{\tau} I_m & 0 & \cdots & 0 \\
0 & \frac{N}{\tau} I_m & -\frac{N}{\tau} I_m & \cdots & 0 \\
\vdots & \vdots & \vdots & \ddots & \vdots \\
0 & 0 & 0 & \cdots & -\frac{N}{\tau} I_m
\end{pmatrix},
\end{align} where in the limit $N\to \infty$
\begin{align*}
{\rm D}F_j(t) = \frac{\partial F(\zeta_0,\zeta_1)}{\partial \zeta_j}\Bigr|_{(x^{\gamma}(t),x^{\gamma}(t-\tau))}.
\end{align*}
The phase adjoint equation for system \eqref{eq:Discretised} is therefore of the form \eqref{eq:adjointZ} where $Z(t) = \tp{(\tp{z_0}(t), \tp{z_1}(t), \ldots, \tp{z_N}(t))}$ and the normalisation condition is 
$\sum_{i=0}^N \tp{z_i}(0) \dot{x}^\gamma_i(0)=\omega$.
 Inserting the substitution $z_0(t) = z(t)$ and \begin{equation} z_i(t) =\frac{\tau}{N} \tp{{\rm D}F_1}\left(t + \tau - (i-1)\tau/N\right)z(t+ \tau - (i-1)\tau/N),\label{eq:zsubs}
\end{equation} for $i =1, \ldots, N$, Novi\v{c}enko and Pyragas \cite{Novicenko2012} observe that in the limit $N \to \infty$, $z(t)$ (the phase response function for the DDE \eqref{eq:DDE}) satisfies the adjoint equation
\begin{equation}
\label{eq:DDEadjointZ}
\dot{z}(t) = -\tp{{\rm D}F_0}(t)z(t) -  \tp{{\rm D}F_1}(t+\tau)z(t+\tau) ,
\end{equation} with the normalisation
\begin{equation}\label{eq:phasenormalisation}
\tp{z}(0) \dot{x}^\gamma(0) + \int_{-\tau}^0 \tp{z}(\tau + \zeta) {\rm D}F_1(\tau +\zeta) \dot{x}^\gamma(\zeta){\rm d} \zeta= \omega .
\end{equation}

Similarly, extending the work of Novi\v{c}enko and Pyragas \cite{Novicenko2012}, the amplitude adjoint equation for system \eqref{eq:Discretised} is of the form \eqref{eq:adjointI} where $I(t) = \tp{(\tp{q_0}(t), \tp{q_1}(t), \ldots, \tp{q_N}(t))}$. Observe that for an ODE system, if we write $W(t) = {\rm e}^{-\mu t} I(t)$ then $W(t)$ satisfies \eqref{eq:adjointZ} when $I(t)$ solves \eqref{eq:adjointI}. Writing $W(t) = \tp{(\tp{w_0}(t), \tp{w_1}(t), \ldots, \tp{w_N}(t))}$, making the substitution $w_0(t)=w(t)$ and $w_i(t)$ as in \eqref{eq:zsubs} with $w$ replacing $z$, and then taking the limit $N \to \infty$ we conclude that $w(t) = w_0(t) = {\rm e}^{-\mu t} q_0(t) = {\rm e}^{-\mu t} q(t)$ satisfies \eqref{eq:DDEadjointZ}. Therefore the amplitude response for the DDE \eqref{eq:DDE}, $q(t)$ satisfies the adjoint equation
\begin{equation}
\label{eq:DDEadjointI}
\dot{q}(t) = -(\tp{{\rm D}F_0}(t) - \mu I_m) q(t) - {\rm e}^{-\mu\tau}\tp{{\rm D}F_1}(t+\tau)q(t+\tau).
\end{equation}
The normalisation for \eqref{eq:DDEadjointI} requires the Floquet eigenfunction $\rho(t)$ for the DDE. For the discretised ODE system \eqref{eq:Discretised} the Floquet eigenfunction $g(t) = \tp{(\tp{\rho_0}(t),\tp{\rho_1}(t), \ldots, \tp{\rho_N}(t))}$ satisfies \eqref{eq:FloquetODE}. Then $h(t) = {\rm e}^{\mu t} g(t)$ satisfies $\dot{h}= Jh$ (the linearised equation for deviations $h$ from limit cycle). In the limit $N\to \infty$ this has solution $h_i(t) = h_0(t-i\tau/N)$ where
\begin{equation}
\dot{h}_0(t) = {\rm D}F_0(t) h_0(t) + {\rm D}F_1(t) h_0(t-\tau).
\end{equation} Therefore $\rho_0(t) = \rho(t)$ satisfies 
\begin{equation}\label{eq:DDEeigenfunctions}
\dot{\rho}(t) = ({\rm D}F_0(t)-\mu I_m) \rho(t) + {\rm e}^{-\mu\tau}{\rm D}F_1(t) \rho(t-\tau),
\end{equation} and we make the normalisation choice $\max_{t} |\rho(t)|=1$. We also note that $\dot{x}^\gamma(t)$ is the eigenfunction of the linearised system with $\mu=0$. 

Since for $i=1, \ldots, N$, $q_i(0)= w_i(0)$ and $\rho_i(0)= h_i(0)$ we see that 
\begin{align*}
\sum_{i=0}^N \tp{q_i}(0)\rho_i(0) = \tp{q}(0) \rho(0) + \sum_{i=1}^N  \tp{w_i}(0)h_i(0).
\end{align*} Taking the limit $N\to \infty$ the normalisation condition for \eqref{eq:DDEadjointI} is given by 
\footnote{
We note that the linear equations \eqref{eq:DDEadjointZ}, \eqref{eq:DDEadjointI} and \eqref{eq:DDEeigenfunctions} derived for the DDE \eqref{eq:DDE} agree with those in \cite{Kotani2020}, however our normalisation for the amplitude response \eqref{eq:normalisation} differs from that in \cite{Kotani2020} by the factor of ${\rm e}^{-\mu \tau}$ multiplying the integral. 
}
\begin{equation}
\tp{q}(0) \rho(0) +{\rm e}^{-\mu \tau}  \int_{-\tau}^0 \tp{q}(\tau+ \zeta) {\rm D}F_1(\tau +\zeta) \rho(\zeta) {\rm d} \zeta =1 . \label{eq:normalisation}
\end{equation}

We now turn to the problem of how to efficiently numerically calculate the quantities $z(t)$, $q(t)$, $\rho(t)$ and $\mu$ as well as the limit cycle. Improvements in accuracy and computational speed can be made over numerical DDE solvers by using the method of harmonic balance (also known as the Fourier–Galerkin method). The harmonic balance method is widely used for analysing the periodic solutions of nonlinear differential equations, often in a mechanical setting, see e.g., \cite{Detroux2015}. The periodic solution in $\mathbb{R}^m$ is expressed as a Fourier series truncated at $M$ Fourier modes. To determine the $m(2M+1)$ unknown Fourier coefficients, an algebraic zero problem of size $m(2M+1)$ is formulated by considering the solution sampled at $2M+1$ time points. The method can also be used for determining periodic solutions of delay differential equations, see e.g., \cite{Sun2023}. Simmendinger \textit{et al}. \cite{Simmendinger1999} also use Fourier expansions to approximately calculate Floquet exponents and eigenfunctions for delay differential equations.

For the DDE \eqref{eq:DDE} first consider the $T$-periodic orbit $x^\gamma(t) \in \mathbb{R}^m$ with truncated Fourier series representation
\begin{equation}
x^\gamma(t) = \sum_{p=-M}^M a_p {\rm e}^{{\rm i} \omega_p t}, \quad \omega_p = \frac{2\pi p}{T} ,
\end{equation} where $a_p \in \mathbb{C}^m$ with $a_{-p} = a_p^*$ and ${}^*$ denotes complex conjugation. There are $m(2M+1)+1$ real unknowns to determine, namely the Fourier coefficients $a_0, a_1, \ldots, a_M$ and the period $T$. Sample $x^\gamma(t)$ at $2M+1$ time instants $t_n = nT/(2M+1)$, $n=-M, \ldots, 0, \ldots, M$ and introduce the notation
\begin{align}
X &= \tp{\begin{bmatrix}
x^{\gamma\mathsf{T}}(t_{-M}) & \cdots & x^{\gamma\mathsf{T}}(t_{0}) & \cdots & x^{\gamma\mathsf{T}}(t_{M})
\end{bmatrix}},\\
A &= \tp{\begin{bmatrix}
\tp{a_{-M}} & \cdots & \tp{a_0} & \cdots & \tp{a_M}
\end{bmatrix}}.
\end{align} Then $X= (S \otimes I_m) A$ where $S$ is the symmetric Vandermonde matrix with 
\begin{equation*} [S]_{np} =  {\rm e}^{2\pi {\rm i} np /(2M+1)},\quad n, p \in \{ -M , \ldots, 0, \ldots M\},\end{equation*} 
and $\otimes$ is the Kronecker product. Sampling \eqref{eq:DDE} on the limit cycle results in the system of $m(2M+1)$ nonlinear algebraic equations for the components of $X$
\begin{equation} \label{eq:zerop1}
((SL(0)S^{-1}) \otimes I_m)X - F(X, ((S\Gamma S^{-1}) \otimes I_m)X)=0,
\end{equation} where $[L(0)]_{np} = \delta_{np} {\rm i} \omega_p$, $[\Gamma]_{np} = \delta_{np} {\rm e}^{-{\rm i} \omega_p \tau}$ and we note that $S^{-1}= S^*/(2M+1)$. To fix an origin of the one parameter family of periodic orbits we set the $n$th component of $x^\gamma$ to have a vanishing derivative at the origin, or equivalently: 
\begin{equation} \label{eq:zerop2}
\tp{e_n}\{ \tp{(-M, \ldots, 0, \ldots, M)}(S^{-1}\otimes I_m) X\}=0,
\end{equation}where $e_n \in \mathbb{R}^m$ is a canonical vector for the $n$th direction. The zero problem given by \eqref{eq:zerop1}-\eqref{eq:zerop2} can then be solved numerically for $X$, for example using an optimisation approach based upon the Levenberg-Marquardt algorithm.

To determine stability we consider again \eqref{eq:DDEeigenfunctions} describing a small $T$-periodic perturbation $\rho(t)$ such that $x(t) = x^\gamma(t) + {\rm e}^{\mu t} \rho(t)$ for $\mu \neq 0$. Sampling \eqref{eq:DDEeigenfunctions} at times $t_n$ and using a truncated Fourier series representation of $\rho(t)$ as $\rho(t) = \sum_{p=-M}^M b_p {\rm e}^{{\rm i} \omega_p t}$ yields the linear algebraic system
\begin{equation}
\label{eq:sampstab} \mathcal{M}(\mu)R=0,
\end{equation}
 for $R = \tp{\begin{bmatrix}\tp{\rho}(t_{-M}) & \cdots & \tp{\rho}(t_{0}) & \cdots & \tp{\rho}(t_{M})\end{bmatrix}}$ and 
\begin{equation*}
\label{eq:Mmu}
\mathcal{M}(\mu) = (S L(\mu) S^{-1}) \otimes I_m - J_0 - {\rm e}^{-\mu \tau} J_1((S \Gamma S^{-1}) \otimes I_m) ,
\end{equation*} where \begin{align*} J_i &= \mbox{Blockdiag}\left({\rm D}F_i(t_{-M}), \ldots, {\rm D}F_i(t_{0}), \ldots, {\rm D}F_i(t_{M}) \right),\end{align*} and $[L(\mu)]_{np} = \delta_{np} (\mu + {\rm i} \omega_p)$ and also $R = (S\otimes I_m)B$ for $B =\tp{\begin{bmatrix}\tp{b_{-M}} & \cdots & \tp{b_0} & \cdots & \tp{b_M}\end{bmatrix}}$. For nontrivial solutions of \eqref{eq:sampstab} we require that $\mathcal{E}(\mu)=\det(\mathcal{M}(\mu))= 0$. The approximated Floquet exponents $\mu_j$ are solutions of $\mathcal{E}(\mu)=0$ with corresponding sampled eigenfunctions $R_j$. For stable limit cycles $\mbox{Re}(\mu_j)<0$ for all nonzero Floquet exponents. Here we consider the case where the largest nontrivial Floquet exponent $\mu_1 := \mu$ is real. 

The phase and amplitude responses can be computed as appropriately normalised solutions of the adjoint equation \eqref{eq:DDEadjointI} (with $\mu=0$ giving phase response). Sampling \eqref{eq:DDEadjointI} at times $t_n$ and using a truncated Fourier series representation of $q(t)$ as $q(t) = \sum_{p=-M}^M c_p {\rm e}^{{\rm i} \omega_p t}$ yields another linear algebraic system
\begin{align}
\label{eq:HBadjoint}
\begin{split}
\tp{Q}\Bigl[ (SL(-\mu)S^{-1})& \otimes I_m   + J_0 +\\ &{\rm e}^{-\mu \tau} ((S\Gamma^* S^{-1}) \otimes I_m)\widetilde{J}_1\Bigr]=0 ,
\end{split}
\end{align}
where 
\begin{align*}
Q &= \tp{\begin{bmatrix}\tp{q}(t_{-M}) & \cdots & \tp{q}(t_{0}) & \cdots & \tp{q}(t_{M})\end{bmatrix}} ,\\
\widetilde{J}_1 &= \mbox{Blockdiag}\left(\widetilde{{\rm D}F}_1(t_{-M}), \ldots, \widetilde{{\rm D}F}_1(t_{0}), \ldots, \widetilde{{\rm D}F}_1(t_{M}) \right) ,\\ &
\widetilde{{\rm D}F}_1(t) = \frac{\partial F(\zeta_0,\zeta_1)}{\partial \zeta_1}\Bigr|_{(x^{\gamma}(t+ \tau),x^{\gamma}(t))} ,
\end{align*}
and the relationship to the Fourier coefficients is given by $Q = (S\otimes I_m)C$ for $C =\tp{\begin{bmatrix}\tp{c_{-M}} & \cdots & \tp{c_0} & \cdots & \tp{c_M}\end{bmatrix}}$. The systems of linear equations \eqref{eq:sampstab} and \eqref{eq:HBadjoint} can be solved faster than using DDE solvers to determine the eigenfunctions and response functions, also recalling the normalisations \eqref{eq:phasenormalisation} when $\mu=0$ and $\max_{t} |\rho(t)|=1$ and \eqref{eq:normalisation} otherwise. 
By comparison with analytically tractable examples (see next), we find that this methodology can give very accurate results.

\begin{figure}[h]
\begin{center}
\includegraphics[width = 0.5\textwidth]{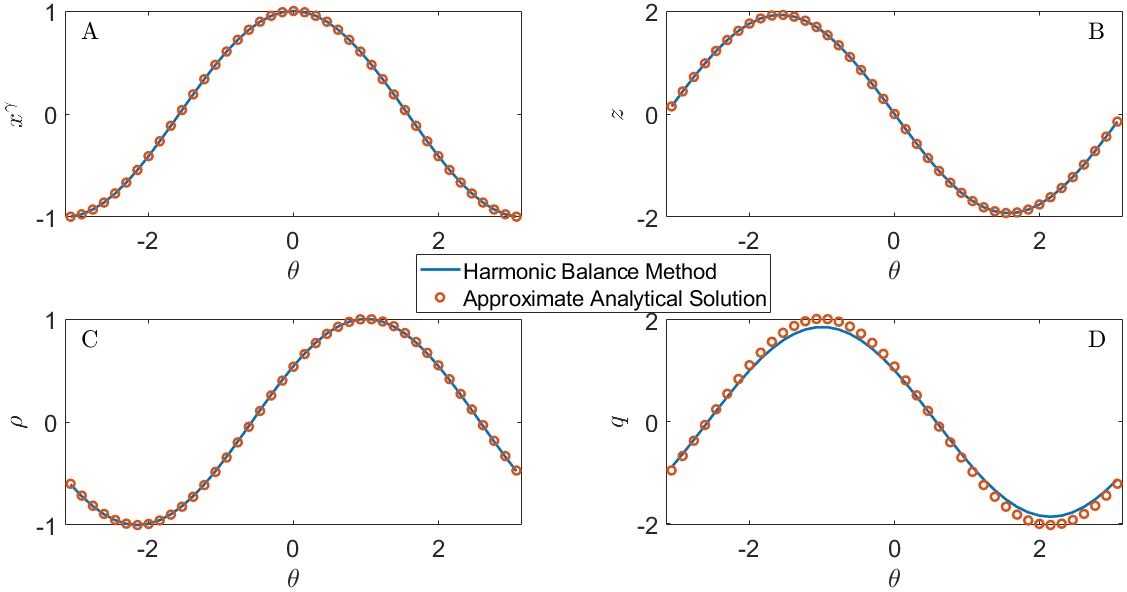}
\caption{The (A) limit cycle, (B) phase response curve (C) Floquet eigenfunction, and (D) amplitude response curve for model \eqref{eq:analyticalDDE} plotted against time.
The analytical curves from \cite{Kotani2020} are shown with red circles and numerical results using the harmonic balance method are shown with blue lines. Here $\delta=0.05$ and $M=20$.} 
\label{fig:HarmvsPaper}
\end{center}
\end{figure}

The simple DDE model 
\begin{equation}
\label{eq:analyticalDDE}
\frac{{\rm d} x(t)}{{\rm d} t} = -x(t-\pi/2) + \delta x(t)\left( 1-x(t)^2 - x(t-\pi/2)^2 \right), 
\end{equation} has a limit cycle $x^\gamma(t) = \cos(t)$. For small values of $\delta$ the forms for $\mu$, $\rho(t)$, $z(t)$, and $q(t)$ can be found analytically \cite{Kotani2020}. For $\epsilon=0.05$ we plot these analytical expressions against the result of the harmonic balance method with truncations at $M=20$ in Figure \ref{fig:HarmvsPaper}. The harmonic balance approach faithfully reproduces the eigenfunction and both phase and amplitude response curves. We also find that harmonic balance improves on the accuracy in addition to the speed of these computations over standard DDE solvers (not shown).

\begin{figure}[h]
\begin{center}
\includegraphics[width =0.5 \textwidth]{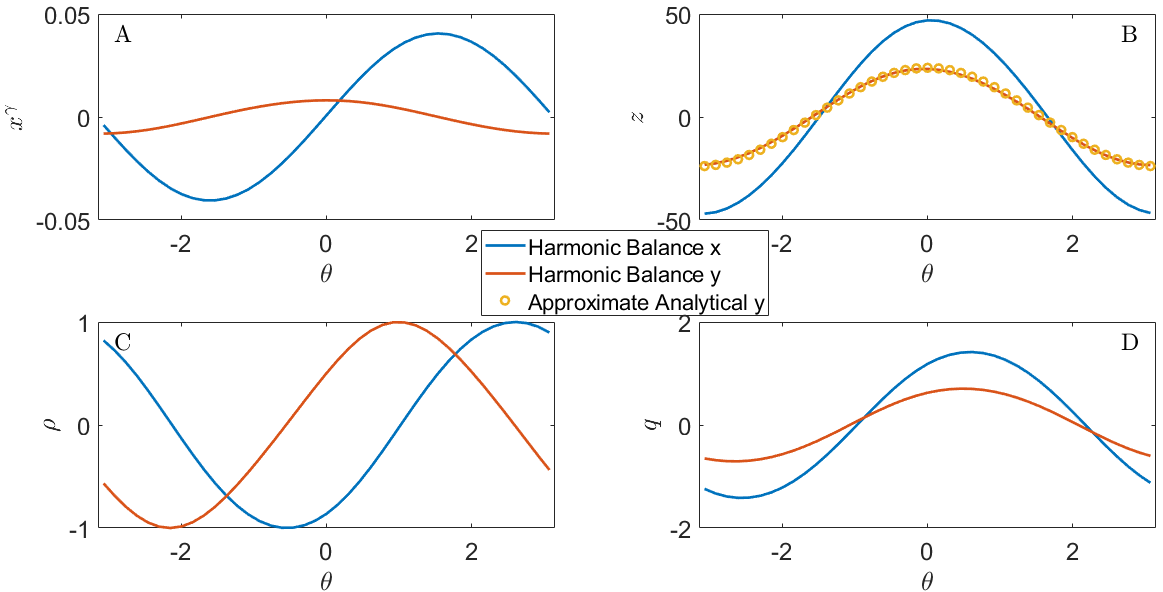}
\caption{ The (A) limit cycle, (B) phase response curve (C) Floquet eigenfunction, and (D) amplitude response curve for model \eqref{eq:CompactCT} plotted against phase. The analytical curve from \cite{Yamaguchi2011} is shown with yellow circles and the $x$ and $y$ components from the harmonic balance method are shown with a blue and red line respectively. Parameter values are $\alpha = -0.039$, $\beta = -0.4$, $\gamma = -2.0$, $\delta = -10.0$ and $\tau = 8.0$ and $M=20$.} 
\label{fig:Cortico}
\end{center}
\end{figure}

As a second (nonscalar) example we consider a simplified model for cortico-thalamic EEG rhythms \cite{Kim2007},
\begin{align}
\frac{{\rm d} x(t)}{{\rm d} t} &= y(t) ,\nonumber\\
\frac{{\rm d} y(t)}{{\rm d} t} &= \gamma y(t) + \alpha x(t) + \beta x(t-\tau) + \delta x(t)^3. \label{eq:CompactCT}
\end{align}
For $\alpha = -0.039$, $\beta = -0.4$, $\gamma = -2.0$, $\delta = -10.0$ and $\tau = 8.0$,  centre manifold reduction can be used to analytically approximate the $y$ component of the phase response \cite{Yamaguchi2011}. Kotani \textit{et al}. \cite{Kotani2012} match this approximation to the PRC calculated by numerically integrating the adjoint equation and by direct perturbation. Using the harmonic balance approach we calculate the limit cycle and the phase response to perturbations in both the $x$ and $y$ components as in Figure \ref{fig:Cortico}(A) and (B) respectively, showing that we also obtain agreement with the analytical expression from \cite{Yamaguchi2011}. We further use harmonic balance to determine the Floquet exponent as $\mu = -0.00296$ and to compute the corresponding eigenfunction and amplitude response function as shown in Figure \ref{fig:Cortico}(C) and (D) respectively. In this case, since the Floquet exponent is so close to zero one would expect only a very slow decay to the limit cycle, which would require long computation times using a direct numerical method based upon time evolution.  This issue does not arise when using the harmonic balance approach, whose accuracy can be further improved with increasingly larger choices for  the truncation parameter $M$.

For simplicity we have restricted to systems with a single time delay, however it is straightforward to extend the approaches given here to systems with multiple time delays. Extensions to consider systems where the Floquet exponent with largest real part is complex are natural and allow for the treatment for a wider range of DDEs. Furthermore, higher order correction terms to the phase and isostable response curves can be computed for ordinary differential equation systems with limit cycle solutions \cite{Wilson2018, Wilson2019, Wilson2019b, Wilson2020b} by solving inhomogeneous linear equations. Analogues of these equations for delay differential equations can also be computed and solved using harmonic balance to obtain more accurate equations for the dynamics of the phase and amplitude of solutions of delay systems. This allows for the exploration of synchronisation and phase-locked state dynamics in networks of delay systems, with and without additional interaction delays, extending previous work \cite{Nicks2024} to delay systems. These extensions will be reported on at length elsewhere.

\bibliographystyle{apsrev4-1}
\bibliography{AmpResponsePRL}% Produces the bibliography via BibTeX.

%merlin.mbs apsrev4-1.bst 2010-07-25 4.21a (PWD, AO, DPC) hacked
%Control: key (0)
%Control: author (72) initials jnrlst
%Control: editor formatted (1) identically to author
%Control: production of article title (-1) disabled
%Control: page (0) single
%Control: year (1) truncated
%Control: production of eprint (0) enabled
\begin{thebibliography}{22}%
\makeatletter
\providecommand \@ifxundefined [1]{%
 \@ifx{#1\undefined}
}%
\providecommand \@ifnum [1]{%
 \ifnum #1\expandafter \@firstoftwo
 \else \expandafter \@secondoftwo
 \fi
}%
\providecommand \@ifx [1]{%
 \ifx #1\expandafter \@firstoftwo
 \else \expandafter \@secondoftwo
 \fi
}%
\providecommand \natexlab [1]{#1}%
\providecommand \enquote  [1]{``#1''}%
\providecommand \bibnamefont  [1]{#1}%
\providecommand \bibfnamefont [1]{#1}%
\providecommand \citenamefont [1]{#1}%
\providecommand \href@noop [0]{\@secondoftwo}%
\providecommand \href [0]{\begingroup \@sanitize@url \@href}%
\providecommand \@href[1]{\@@startlink{#1}\@@href}%
\providecommand \@@href[1]{\endgroup#1\@@endlink}%
\providecommand \@sanitize@url [0]{\catcode `\\12\catcode `\$12\catcode
  `\&12\catcode `\#12\catcode `\^12\catcode `\_12\catcode `\%12\relax}%
\providecommand \@@startlink[1]{}%
\providecommand \@@endlink[0]{}%
\providecommand \url  [0]{\begingroup\@sanitize@url \@url }%
\providecommand \@url [1]{\endgroup\@href {#1}{\urlprefix }}%
\providecommand \urlprefix  [0]{URL }%
\providecommand \Eprint [0]{\href }%
\providecommand \doibase [0]{http://dx.doi.org/}%
\providecommand \selectlanguage [0]{\@gobble}%
\providecommand \bibinfo  [0]{\@secondoftwo}%
\providecommand \bibfield  [0]{\@secondoftwo}%
\providecommand \translation [1]{[#1]}%
\providecommand \BibitemOpen [0]{}%
\providecommand \bibitemStop [0]{}%
\providecommand \bibitemNoStop [0]{.\EOS\space}%
\providecommand \EOS [0]{\spacefactor3000\relax}%
\providecommand \BibitemShut  [1]{\csname bibitem#1\endcsname}%
\let\auto@bib@innerbib\@empty
%</preamble>
\bibitem [{\citenamefont {Erneux}(2009)}]{Erneux2009}%
  \BibitemOpen
  \bibfield  {author} {\bibinfo {author} {\bibfnamefont {T.}~\bibnamefont
  {Erneux}},\ }\href@noop {} {\emph {\bibinfo {title} {Applied Delay
  Differential Equations}}},\ Surveys and Tutorials in the Applied Mathematical
  Sciences\ (\bibinfo  {publisher} {Springer},\ \bibinfo {year}
  {2009})\BibitemShut {NoStop}%
\bibitem [{\citenamefont {Stepan}(1989)}]{Stepan1989}%
  \BibitemOpen
  \bibfield  {author} {\bibinfo {author} {\bibfnamefont {G.}~\bibnamefont
  {Stepan}},\ }\href@noop {} {\emph {\bibinfo {title} {Retarded Dynamical
  Systems: Stability and Characteristic Functions}}}\ (\bibinfo  {publisher}
  {Longman Higher Education},\ \bibinfo {year} {1989})\BibitemShut {NoStop}%
\bibitem [{\citenamefont {Campbell}(2007)}]{Campbell2007}%
  \BibitemOpen
  \bibfield  {author} {\bibinfo {author} {\bibfnamefont {S.~A.}\ \bibnamefont
  {Campbell}},\ }\enquote {\bibinfo {title} {Time delays in neural systems},}\
  in\ \href {\doibase 10.1007/978-3-540-71512-2_2} {\emph {\bibinfo {booktitle}
  {Handbook of Brain Connectivity}}},\ \bibinfo {editor} {edited by\ \bibinfo
  {editor} {\bibfnamefont {V.~K.}\ \bibnamefont {Jirsa}}\ and\ \bibinfo
  {editor} {\bibfnamefont {A.~R.}\ \bibnamefont {McIntosh}}}\ (\bibinfo
  {publisher} {Springer Berlin Heidelberg},\ \bibinfo {address} {Berlin,
  Heidelberg},\ \bibinfo {year} {2007})\ pp.\ \bibinfo {pages}
  {65--90}\BibitemShut {NoStop}%
\bibitem [{\citenamefont {Guillamon}\ and\ \citenamefont
  {Huguet}(2009)}]{Guillamon2009}%
  \BibitemOpen
  \bibfield  {author} {\bibinfo {author} {\bibfnamefont {A.}~\bibnamefont
  {Guillamon}}\ and\ \bibinfo {author} {\bibfnamefont {G.}~\bibnamefont
  {Huguet}},\ }\href {\doibase 10.1137/080737666} {\bibfield  {journal}
  {\bibinfo  {journal} {SIAM Journal on Applied Dynamical Systems}\ }\textbf
  {\bibinfo {volume} {8}},\ \bibinfo {pages} {1005} (\bibinfo {year}
  {2009})}\BibitemShut {NoStop}%
\bibitem [{\citenamefont {Mauroy}\ \emph {et~al.}(2013)\citenamefont {Mauroy},
  \citenamefont {Mezi{\'c}},\ and\ \citenamefont {Moehlis}}]{Mauroy2013}%
  \BibitemOpen
  \bibfield  {author} {\bibinfo {author} {\bibfnamefont {A.}~\bibnamefont
  {Mauroy}}, \bibinfo {author} {\bibfnamefont {I.}~\bibnamefont {Mezi{\'c}}}, \
  and\ \bibinfo {author} {\bibfnamefont {J.}~\bibnamefont {Moehlis}},\ }\href
  {\doibase https://doi.org/10.1016/j.physd.2013.06.004} {\bibfield  {journal}
  {\bibinfo  {journal} {Physica D}\ }\textbf {\bibinfo {volume} {261}},\
  \bibinfo {pages} {19} (\bibinfo {year} {2013})}\BibitemShut {NoStop}%
\bibitem [{\citenamefont {Wilson}\ and\ \citenamefont
  {Moehlis}(2016)}]{Wilson2016}%
  \BibitemOpen
  \bibfield  {author} {\bibinfo {author} {\bibfnamefont {D.}~\bibnamefont
  {Wilson}}\ and\ \bibinfo {author} {\bibfnamefont {J.}~\bibnamefont
  {Moehlis}},\ }\href {\doibase 10.1103/PhysRevE.94.052213} {\bibfield
  {journal} {\bibinfo  {journal} {Physical Review E}\ }\textbf {\bibinfo
  {volume} {94}},\ \bibinfo {pages} {052213} (\bibinfo {year}
  {2016})}\BibitemShut {NoStop}%
\bibitem [{\citenamefont {Ermentrout}\ \emph {et~al.}(2019)\citenamefont
  {Ermentrout}, \citenamefont {Park},\ and\ \citenamefont
  {Wilson}}]{Ermentrout2019}%
  \BibitemOpen
  \bibfield  {author} {\bibinfo {author} {\bibfnamefont {B.}~\bibnamefont
  {Ermentrout}}, \bibinfo {author} {\bibfnamefont {Y.}~\bibnamefont {Park}}, \
  and\ \bibinfo {author} {\bibfnamefont {D.}~\bibnamefont {Wilson}},\
  }\href@noop {} {\bibfield  {journal} {\bibinfo  {journal} {Philosophical
  Transactions of the Royal Society A}\ }\textbf {\bibinfo {volume} {377}},\
  \bibinfo {pages} {20190092} (\bibinfo {year} {2019})}\BibitemShut {NoStop}%
\bibitem [{\citenamefont {Nicks}\ \emph {et~al.}(2024)\citenamefont {Nicks},
  \citenamefont {Allen},\ and\ \citenamefont {Coombes}}]{Nicks2024}%
  \BibitemOpen
  \bibfield  {author} {\bibinfo {author} {\bibfnamefont {R.}~\bibnamefont
  {Nicks}}, \bibinfo {author} {\bibfnamefont {R.}~\bibnamefont {Allen}}, \ and\
  \bibinfo {author} {\bibfnamefont {S.}~\bibnamefont {Coombes}},\ }\href
  {\doibase 10.1063/5.0179430} {\bibfield  {journal} {\bibinfo  {journal}
  {Chaos}\ }\textbf {\bibinfo {volume} {34}},\ \bibinfo {pages} {013141}
  (\bibinfo {year} {2024})}\BibitemShut {NoStop}%
\bibitem [{\citenamefont {Novi\v{c}enko}\ and\ \citenamefont
  {Pyragas}(2012)}]{Novicenko2012}%
  \BibitemOpen
  \bibfield  {author} {\bibinfo {author} {\bibfnamefont {V.}~\bibnamefont
  {Novi\v{c}enko}}\ and\ \bibinfo {author} {\bibfnamefont {K.}~\bibnamefont
  {Pyragas}},\ }\href@noop {} {\bibfield  {journal} {\bibinfo  {journal}
  {Physica D}\ }\textbf {\bibinfo {volume} {241}},\ \bibinfo {pages} {1090}
  (\bibinfo {year} {2012})}\BibitemShut {NoStop}%
\bibitem [{\citenamefont {Kotani}\ \emph {et~al.}(2012)\citenamefont {Kotani},
  \citenamefont {Yamaguchi}, \citenamefont {Ogawa}, \citenamefont {Jimbo},
  \citenamefont {Nakao},\ and\ \citenamefont {Ermentrout}}]{Kotani2012}%
  \BibitemOpen
  \bibfield  {author} {\bibinfo {author} {\bibfnamefont {K.}~\bibnamefont
  {Kotani}}, \bibinfo {author} {\bibfnamefont {I.}~\bibnamefont {Yamaguchi}},
  \bibinfo {author} {\bibfnamefont {Y.}~\bibnamefont {Ogawa}}, \bibinfo
  {author} {\bibfnamefont {Y.}~\bibnamefont {Jimbo}}, \bibinfo {author}
  {\bibfnamefont {H.}~\bibnamefont {Nakao}}, \ and\ \bibinfo {author}
  {\bibfnamefont {G.~B.}\ \bibnamefont {Ermentrout}},\ }\href@noop {}
  {\bibfield  {journal} {\bibinfo  {journal} {Physical Review Letters}\
  }\textbf {\bibinfo {volume} {109}},\ \bibinfo {pages} {044101} (\bibinfo
  {year} {2012})}\BibitemShut {NoStop}%
\bibitem [{\citenamefont {Kotani}\ \emph {et~al.}(2020)\citenamefont {Kotani},
  \citenamefont {Ogawa}, \citenamefont {Shirasaka}, \citenamefont {Akao},
  \citenamefont {Jimbo},\ and\ \citenamefont {Nakao}}]{Kotani2020}%
  \BibitemOpen
  \bibfield  {author} {\bibinfo {author} {\bibfnamefont {K.}~\bibnamefont
  {Kotani}}, \bibinfo {author} {\bibfnamefont {Y.}~\bibnamefont {Ogawa}},
  \bibinfo {author} {\bibfnamefont {S.}~\bibnamefont {Shirasaka}}, \bibinfo
  {author} {\bibfnamefont {A.}~\bibnamefont {Akao}}, \bibinfo {author}
  {\bibfnamefont {Y.}~\bibnamefont {Jimbo}}, \ and\ \bibinfo {author}
  {\bibfnamefont {H.}~\bibnamefont {Nakao}},\ }\href {\doibase
  10.1103/PhysRevResearch.2.033106} {\bibfield  {journal} {\bibinfo  {journal}
  {Physical Review Research}\ }\textbf {\bibinfo {volume} {2}},\ \bibinfo
  {pages} {033106} (\bibinfo {year} {2020})}\BibitemShut {NoStop}%
\bibitem [{\citenamefont {Wilson}\ and\ \citenamefont
  {Ermentrout}(2018)}]{Wilson2018}%
  \BibitemOpen
  \bibfield  {author} {\bibinfo {author} {\bibfnamefont {D.}~\bibnamefont
  {Wilson}}\ and\ \bibinfo {author} {\bibfnamefont {B.}~\bibnamefont
  {Ermentrout}},\ }\href {\doibase 10.1007/s00285-017-1141-6} {\bibfield
  {journal} {\bibinfo  {journal} {Journal of Mathematical Biology}\ }\textbf
  {\bibinfo {volume} {76}},\ \bibinfo {pages} {37} (\bibinfo {year}
  {2018})}\BibitemShut {NoStop}%
\bibitem [{\citenamefont {Park}\ and\ \citenamefont {Wilson}(2021)}]{Park2021}%
  \BibitemOpen
  \bibfield  {author} {\bibinfo {author} {\bibfnamefont {Y.}~\bibnamefont
  {Park}}\ and\ \bibinfo {author} {\bibfnamefont {D.~D.}\ \bibnamefont
  {Wilson}},\ }\href {\doibase 10.1137/20M1371208} {\bibfield  {journal}
  {\bibinfo  {journal} {SIAM Journal on Applied Dynamical Systems}\ }\textbf
  {\bibinfo {volume} {20}},\ \bibinfo {pages} {1464} (\bibinfo {year}
  {2021})}\BibitemShut {NoStop}%
\bibitem [{\citenamefont {Wilson}(2020)}]{Wilson2020b}%
  \BibitemOpen
  \bibfield  {author} {\bibinfo {author} {\bibfnamefont {D.}~\bibnamefont
  {Wilson}},\ }\href {\doibase 10.1103/PhysRevE.101.022220} {\bibfield
  {journal} {\bibinfo  {journal} {Physical Review E}\ }\textbf {\bibinfo
  {volume} {101}},\ \bibinfo {pages} {022220} (\bibinfo {year}
  {2020})}\BibitemShut {NoStop}%
\bibitem [{Note1()}]{Note1}%
  \BibitemOpen
  \bibinfo {note} {We note that the linear equations \protect \textup {\hbox
  {\mathsurround \z@ \protect \normalfont (\ignorespaces \ref
  {eq:DDEadjointZ}\unskip \@@italiccorr )}}, \protect \textup {\hbox
  {\mathsurround \z@ \protect \normalfont (\ignorespaces \ref
  {eq:DDEadjointI}\unskip \@@italiccorr )}} and \protect \textup {\hbox
  {\mathsurround \z@ \protect \normalfont (\ignorespaces \ref
  {eq:DDEeigenfunctions}\unskip \@@italiccorr )}} derived for the DDE \protect
  \textup {\hbox {\mathsurround \z@ \protect \normalfont (\ignorespaces \ref
  {eq:DDE}\unskip \@@italiccorr )}} agree with those in \cite {Kotani2020},
  however our normalisation for the amplitude response \protect \textup {\hbox
  {\mathsurround \z@ \protect \normalfont (\ignorespaces \ref
  {eq:normalisation}\unskip \@@italiccorr )}} differs from that in \cite
  {Kotani2020} by the factor of ${\protect \rm e}^{-\mu \tau }$ multiplying the
  integral.}\BibitemShut {Stop}%
\bibitem [{\citenamefont {Detroux}\ \emph {et~al.}(2015)\citenamefont
  {Detroux}, \citenamefont {Renson}, \citenamefont {Masset},\ and\
  \citenamefont {Kerschen}}]{Detroux2015}%
  \BibitemOpen
  \bibfield  {author} {\bibinfo {author} {\bibfnamefont {T.}~\bibnamefont
  {Detroux}}, \bibinfo {author} {\bibfnamefont {L.}~\bibnamefont {Renson}},
  \bibinfo {author} {\bibfnamefont {L.}~\bibnamefont {Masset}}, \ and\ \bibinfo
  {author} {\bibfnamefont {G.}~\bibnamefont {Kerschen}},\ }\href {\doibase
  https://doi.org/10.1016/j.cma.2015.07.017} {\bibfield  {journal} {\bibinfo
  {journal} {Computer Methods in Applied Mechanics and Engineering}\ }\textbf
  {\bibinfo {volume} {296}},\ \bibinfo {pages} {18} (\bibinfo {year}
  {2015})}\BibitemShut {NoStop}%
\bibitem [{\citenamefont {Sun}\ \emph {et~al.}(2023)\citenamefont {Sun},
  \citenamefont {Zhao}, \citenamefont {Yu}, \citenamefont {Huang},
  \citenamefont {Feng},\ and\ \citenamefont {Zhou}}]{Sun2023}%
  \BibitemOpen
  \bibfield  {author} {\bibinfo {author} {\bibfnamefont {P.}~\bibnamefont
  {Sun}}, \bibinfo {author} {\bibfnamefont {X.}~\bibnamefont {Zhao}}, \bibinfo
  {author} {\bibfnamefont {X.}~\bibnamefont {Yu}}, \bibinfo {author}
  {\bibfnamefont {Q.}~\bibnamefont {Huang}}, \bibinfo {author} {\bibfnamefont
  {Z.}~\bibnamefont {Feng}}, \ and\ \bibinfo {author} {\bibfnamefont
  {J.}~\bibnamefont {Zhou}},\ }\href {\doibase
  https://doi.org/10.1016/j.apm.2023.02.018} {\bibfield  {journal} {\bibinfo
  {journal} {Applied Mathematical Modelling}\ }\textbf {\bibinfo {volume}
  {118}},\ \bibinfo {pages} {818} (\bibinfo {year} {2023})}\BibitemShut
  {NoStop}%
\bibitem [{\citenamefont {Simmendinger}\ \emph {et~al.}(1999)\citenamefont
  {Simmendinger}, \citenamefont {Wunderlin},\ and\ \citenamefont
  {Pelster}}]{Simmendinger1999}%
  \BibitemOpen
  \bibfield  {author} {\bibinfo {author} {\bibfnamefont {C.}~\bibnamefont
  {Simmendinger}}, \bibinfo {author} {\bibfnamefont {A.}~\bibnamefont
  {Wunderlin}}, \ and\ \bibinfo {author} {\bibfnamefont {A.}~\bibnamefont
  {Pelster}},\ }\href {\doibase 10.1103/PhysRevE.59.5344} {\bibfield  {journal}
  {\bibinfo  {journal} {Physical Review E}\ }\textbf {\bibinfo {volume} {59}},\
  \bibinfo {pages} {5344} (\bibinfo {year} {1999})}\BibitemShut {NoStop}%
\bibitem [{\citenamefont {Kim}\ and\ \citenamefont {Robinson}(2007)}]{Kim2007}%
  \BibitemOpen
  \bibfield  {author} {\bibinfo {author} {\bibfnamefont {J.~W.}\ \bibnamefont
  {Kim}}\ and\ \bibinfo {author} {\bibfnamefont {P.~A.}\ \bibnamefont
  {Robinson}},\ }\href {\doibase 10.1103/PhysRevE.75.031907} {\bibfield
  {journal} {\bibinfo  {journal} {Physical Review E}\ }\textbf {\bibinfo
  {volume} {75}},\ \bibinfo {pages} {031907} (\bibinfo {year}
  {2007})}\BibitemShut {NoStop}%
\bibitem [{\citenamefont {Yamaguchi}\ \emph {et~al.}(2011)\citenamefont
  {Yamaguchi}, \citenamefont {Ogawa}, \citenamefont {Jimbo}, \citenamefont
  {Nakao},\ and\ \citenamefont {Kotani}}]{Yamaguchi2011}%
  \BibitemOpen
  \bibfield  {author} {\bibinfo {author} {\bibfnamefont {I.}~\bibnamefont
  {Yamaguchi}}, \bibinfo {author} {\bibfnamefont {Y.}~\bibnamefont {Ogawa}},
  \bibinfo {author} {\bibfnamefont {Y.}~\bibnamefont {Jimbo}}, \bibinfo
  {author} {\bibfnamefont {H.}~\bibnamefont {Nakao}}, \ and\ \bibinfo {author}
  {\bibfnamefont {K.}~\bibnamefont {Kotani}},\ }\href {\doibase
  10.1371/journal.pone.0026497} {\bibfield  {journal} {\bibinfo  {journal}
  {PLoS ONE}\ }\textbf {\bibinfo {volume} {6}},\ \bibinfo {pages} {1} (\bibinfo
  {year} {2011})}\BibitemShut {NoStop}%
\bibitem [{\citenamefont {Wilson}(2019)}]{Wilson2019}%
  \BibitemOpen
  \bibfield  {author} {\bibinfo {author} {\bibfnamefont {D.}~\bibnamefont
  {Wilson}},\ }\href {\doibase 10.1103/PhysRevE.99.022210} {\bibfield
  {journal} {\bibinfo  {journal} {Physical Review E}\ }\textbf {\bibinfo
  {volume} {99}},\ \bibinfo {pages} {022210} (\bibinfo {year}
  {2019})}\BibitemShut {NoStop}%
\bibitem [{\citenamefont {Wilson}\ and\ \citenamefont
  {Ermentrout}(2019)}]{Wilson2019b}%
  \BibitemOpen
  \bibfield  {author} {\bibinfo {author} {\bibfnamefont {D.}~\bibnamefont
  {Wilson}}\ and\ \bibinfo {author} {\bibfnamefont {B.}~\bibnamefont
  {Ermentrout}},\ }\href {\doibase 10.1137/18M1170558} {\bibfield  {journal}
  {\bibinfo  {journal} {SIAM Review}\ }\textbf {\bibinfo {volume} {61}},\
  \bibinfo {pages} {277} (\bibinfo {year} {2019})}\BibitemShut {NoStop}%
\end{thebibliography}%

\end{document}